\documentclass[a4paper,12pt]{article}
\parskip=\medskipamount
\usepackage{amsmath}
\usepackage{amsfonts}
\usepackage{url}
\usepackage{hyperref}
\usepackage{color}
\newcommand\vartextvisiblespace[1][.5em]{%
  \makebox[#1]{%
    \kern.07em
    \vrule height.3ex
    \hrulefill
    \vrule height.3ex
    \kern.07em
  }
}
\begin{document}
\begin{center}
            On the examples of {E}gyptian fractions in  {L}iber {A}baci\\

                          by\\[2ex]
                  Trond Steihaug and Milo Gardner\\
                       February 7, 2025

\end{center}
\section*{Abstract}
The focus of this note is to formulate the algorithms and give the examples used by Fibonacci in Liber Abaci to expand any fraction into a sum of unit fractions. The description in Liber Abaci is all verbal and the examples are numbers which may lead to different algorithmic descriptions with the same input and results. An additional complication is that the manuscripts that exist are copies of an older manuscript and in the process new errors are introduced. Additional errors may also be introduced in the transcript and the translation. Fibonacci introduces seven categories each with several numerical examples. We give a precise description of the computational procedure using standard mathematical notation.

\section*{ Introduction}

Fibonacci's (or Leonardo of Pisa)  treatment of Egyptian fraction is found in {\em Liber Abaci} written around 1202. The manuscript was transcribed by Boncompagni in 1857 \cite{Boncompagni1857} and a translation to English of section on Egyptian fraction was published by Dunton and Grimm in 1966 \cite{Dunton1966}. The complete Liber Abaci was translated to English with extensive comments by Sigler in 2002 \cite{Sigler2002}. Both of these were based on the transcript by Boncompagni of the manuscript  {\em Conventi Soppressi
C.I.2616, in Biblioteca Nazionale Centrale di Firenze}
\cite{Fibonacci_LA_2616}. A revised version of Liber Abaci from around 1220 was transcribed in 2020 by Giusti \cite{Giusti2020}. This transcript contains variations found in the 19 different manuscript where 9 contains all or most of Liber Abaci. All surviving complete manuscripts seem to derive from a single archetype which again is different from the original autograph. In the paper {\em New Editorial Perspectives on Fibonacci’s Liber Abaci} Giuseppe Germano \cite{Germano2013} writes on the translation by Sigler
\begin{quote}
Unfortunately, however, it has succeeded in inevitably
adding to the already numerous mistakes and problems of its original
a series of over-simplifications and misunderstandings, \ldots  but
also because it exhibits in some cases unreliable, if not, at times, bizarre interpretations,
\ldots
\end{quote}

Dunton and Grimm 1966 \cite{Dunton1966} corrects some of the errors in the transcript by Boncompagni but also introduces new ones. The major difference between the two translations is that in \cite{Sigler2002} the original fractional notation is kept, while in \cite{Dunton1966} the modern notation is used.

The focus of this note is to understand the algorithms used by Fibonacci to expand any fraction into a sum of unit fractions.

\subsection*{Egyptian Fraction}

A fraction written as sum of unit fractions where all the unit fractions are different is an Egyptian fraction. Let $p$ and $q$ be positive integers $p<q$, then an Egyptian fraction will be
\begin{equation}\label{eg:EgyptianFraction}
\frac{p}{q}= \frac{1}{d_1}+\frac{1}{d_2}+\cdots+ \frac{1}{d_\ell}
\end{equation}
where $\ell\geq 2$ and sorted $d_i<d_{i+1},\ i=1,\ldots, \ell-1$.\footnote{In ancient Egyptian mathematics also the fraction $\frac{2}{3}$ is used.} Let $p_1=p$ and $q_1=q$, and define the iteration
\begin{equation}\label{eg:EgyptianFraction_2}
\frac{p_{i+1}}{q_{i+1}}= \frac{p_i}{q_i}-\frac{1}{d_i}.
\end{equation}
Then $d_i$ must satisfy $d_i > \frac{q_i}{p_i}$. Assuming that the fraction $p_{i+1}/q_{i+1}$ is reduced ($p_{i+1}$ and $q_{i+1}$  have no common factors.) The iteration terminates when $p_\ell=1$ with $d_\ell=q_k$
 The Fibonacci greedy algorithm will choose
$d_i=\lceil \frac{q_i}{p_i}\rceil$. Here $ \lceil \cdot\rceil$ is the ceiling function, i.e.
$ \lceil x\rceil$ is the smallest integer greater or equal $x$. The greedy algorithm goes also under the name of
Fibonacci-Sylvester expansion \cite{Mays1987}.

\subsection*{Examples}

Consider the two examples in Table \ref{table:examples}.
For $\frac{7}{15}$ the first expansion is the greedy algorithm while the second
is using Fibonacci category 6.
The fraction $\frac{20}{53}$ is a sum of $\frac{18}{53}$ and $\frac{2}{53}$. For each term, Fibonacci category 3 will give a two term Egyptian fraction, $\frac{18}{53} = \frac{1}{3}+ \frac{1}{159}$ and
$\frac{2}{53}=\frac{1}{27}+\frac{1}{1431}$, and $\frac{20}{53}$ will have a four term expansion. However, a three term expansion can be found by
\[\frac{20}{53}=\frac{1}{4}+(\frac{20}{53}-\frac{1}{4})=
\frac{1}{4}+\frac{1}{4}\;\frac{27}{53}=\frac{1}{4}+\frac{1}{4}(\frac{1}{2}+\frac{1}{106})\]
using the third category. The non-uniqueness will lead to question like what is the shortest possible expansion or the one with the smallest denominators. But this also raises the question on the number of possible expansions in unit fractions. For further discussion see \cite{Steihaug2025a}.

\begin{table}[ht!]
\centering
\begin{tabular}{|c|c|c|c|c|}\hline
                    & $i=1$          & $i=2$          & $i=3$ & $i=4$ \\\hline
  $\frac{p_i}{q_i}$ & $\frac{7}{15}$ & $\frac{2}{15}$ & $\frac{1}{120}$ &\\\hline
  $d_i$             & 3              & 8               & 120 &\\ \hline
  $\frac{p_i}{q_i}$ & $\frac{7}{15}$ & $\frac{13}{60}$ & $\frac{1}{60}$ &\\\hline
  $d_i$             & 4              & 5               & 60 &\\ \hline
  \hline
  $\frac{p_i}{q_i}$ & $\frac{20}{53}$ & $\frac{7}{159}$ & $\frac{2}{3675}$&$\frac{1}{6688653}$ \\\hline
    $d_i$             & 3              & 23               & 1829 &6688653\\ \hline
    $\frac{p_i}{q_i}$ & $\frac{20}{53}$ & $\frac{7}{159}$ & $\frac{10}{1431}$&$\frac{1}{1431}$ \\\hline
    $d_i$             & 3              & 27               & 159 &1431\\ \hline
\end{tabular}
\caption{Examples of unit fraction expansion.}
\label{table:examples}
\end{table}

\section*{Egyptian fractions by  Fibonacci }
The material on unit fractions is in the last part of Chapter 7 of Liber Abaci, first transcribed by Boncompagni in 1857 \cite[p.77-83]{Boncompagni1857} and in 2020 by Giusti \cite[p.131-139]{Giusti2020}. Dunton and Grimm \cite{Dunton1966} gave the first English translation of this section solely based on the transcript by Boncompagni.  Sigler's translation follows the Boncompagni transcripts \cite[p.119-126]{Sigler2002}. Lüneburg \cite[p.81-85]{Luneburg1993} has  discussion of parts of the section on unit fractions.

The error reported in the part on Egyptian fractions of Chapter 7, Sigler \cite[p.620]{Sigler2002} attributes this to Fibonacci but is one instance of an error that can not be attributed to the author \cite[p.lxxxii]{Giusti2020}.

Fibonacci describes the derivation of unit fractions with  seven different techniques in English called {\em distinction} \cite{Sigler2002} or {\em categories} \cite{Dunton1966} or in Latin {\em differentia} \cite{Boncompagni1857}. A further discussion of Fibonacci's methods and the reconstruction of ancient Egyptian mathematics in \cite{Steihaug2025a}.

\subsection*{The first category}
In the first category Fibonacci consider fractions on the form $\frac{p}{q}$ or
$\frac{\frac{p}{q}}{r}$. The first category is divided in three parts called simple, composite, and reversed composite.
\begin{itemize}
\item The simple case is when $p$ divides $q$. Then there exists integer $k>1$ so that $q=kp$ and $\frac{p}{q}=\frac{1}{k}$. Fibonacci gives the examples $\frac{3}{12}=\frac{1}{4}$,
    $\frac{4}{20}=\frac{1}{5}$, and $\frac{5}{100}=\frac{1}{20}$.
\item In the composite case consider $\frac{\frac{p}{q}}{r}$ and $p$ divides $q$, $q=kp$,  then
    \[\frac{\frac{p}{q}}{r}=\frac{\frac{1}{k}}{r}=\frac{1}{kr}.\]
     Fibonacci consider three examples
    \[\frac{\frac{2}{4}}{9}=\frac{\frac{1}{2}}{9}=\frac{1}{18},\
    \frac{\frac{2}{6}}{9}=\frac{\frac{1}{3}}{9}, \text{ and } \frac{\frac{3}{9}}{10}=\frac{\frac{1}{3}}{10}.\]
\item In the reversed composite case consider $\frac{\frac{p}{q}}{r}$ and $p$ divides $r$ then there exits $k>1$ so that $r=kp$, and
    \[\frac{\frac{p}{q}}{r}=\frac{\frac{p}{r}}{q}=\frac{\frac{1}{k}}{q}=\frac{1}{kq}.\]
    Fibonacci consider three examples
    \[\frac{\frac{3}{5}}{9}=\frac{\frac{3}{9}}{5}=\frac{\frac{1}{3}}{5},\
    \frac{\frac{4}{7}}{8}=\frac{\frac{4}{8}}{7}=\frac{\frac{1}{2}}{7}, \text{ and } \frac{\frac{5}{9}}{10}=\frac{\frac{5}{10}}{9}=\frac{\frac{1}{2}}{9}.\]
\end{itemize}
\subsection*{The second category}
The second category is divided in three parts.
\begin{itemize}
\item The fraction $\frac{p}{q}$ is irreducible, but $p$ may be written as a sum of distinct positive numbers that divides $q$. Let $p=s_1+s_2+\cdots+s_\ell$
    where $s_i$ divides $q$, i.e. there exists $k_i$ so that $q=s_ik_i$ for each $i$,
    \[\frac{p}{q}= \frac{s_1+s_2+\cdots+s_\ell}{q}= \frac{1}{k_1}+ \frac{1}{k_2}+\cdots +\frac{1}{k_\ell}.\]
    Sigler  \cite[p.620]{Sigler2002} expresses the second category as
    \[\frac{s_1+s_2}{ks_1s_2}=\frac{1}{s_2k}+\frac{1}{s_1k}\]
    which is correct in general only if $s_1$ and $s_2$ are coprime.

    Fibonacci gives two examples
    \[\frac{5}{6}= \frac{2+3}{6} = \frac{1}{3}+\frac{1}{2} \text{ and }
      \frac{7}{8}= \frac{1+2+4}{8} = \frac{1}{8}+\frac{1}{4}+\frac{1}{2}.\]
    Fibonacci provides tables of these representations for fractions $\frac{p}{q}$ having as denominators $q$ the practical numbers 6, 8, 12, 20, 24, 60, and 100 and $1\leq p \leq q$.  We say
    that $q$ is a practical number  if all smaller positive integers $p$  can be represented as sums of distinct divisors of $q$. Fibonacci calls these tables for {\em Table of separation}. There are some fractions that are not written as a sum of unit fractions in the tables like $\frac{2}{5},\ \frac{3}{5},\ \frac{4}{5},\ \frac{3}{4}$, and $\frac{2}{25}$. In the tables are 2/3, 2/5 and 2/15 which are the same sum of unit fractions as in Rhind Mathematical Papyrus 2/$n$ table. This is further discussed in section Tables of separations.
\item A composite case $\frac{\frac{p}{q}}{r}$ where $p$ may be written as a sum of positive numbers that divides $q$. Let $p=s_1+s_2+\cdots+s_\ell$
    where $s_i$ divides $q$, i.e. there exists $k_i$ so that $q=s_ik_i$ for each $i$ then
    \[\frac{\frac{p}{q}}{r}=
        \frac{\frac{s_1+s_2+\cdots+s_\ell}{q}}{r}=
        \frac{\frac{1}{k_1}}{r}+\frac{\frac{1}{k_2}}{r}+\cdots+\frac{\frac{1}{k_\ell}}{r}.\]
    In the first examples Fibonacci uses $3=1+2$ and $5=1+4$ in the second example.
    \[\frac{\frac{3}{4}}{10}= \frac{\frac{1}{4}}{10}+\frac{\frac{2}{4}}{10}=\frac{1}{20}+\frac{1}{40}\text{ and } \frac{\frac{5}{8}}{9}=\frac{\frac{1}{2}}{9}+\frac{\frac{1}{8}}{9}.\]
\item In the reversed composite case $\frac{\frac{p}{q}}{r}$ where $p$ may be written as a sum of positive numbers that divides $r$. Let $p=s_1+s_2+\cdots+s_\ell$
    where $s_i$ divides $r$, i.e. there exists $k_i$ so that $r=s_ik_i$ for each $i$ then
    \[\frac{\frac{p}{q}}{r}= \frac{\frac{s_1+s_2+\cdots+s_\ell}{r}}{q}=
    \frac{\frac{1}{k_1}}{q}+\frac{\frac{1}{k_2}}{q}+\cdots+\frac{\frac{1}{k_\ell}}{q}.\]

    Three examples
    \begin{eqnarray*}
    \frac{\frac{5}{8}}{10}&=& \frac{\frac{5}{10}}{8}= \frac{\frac{1}{2}}{8},\\ \frac{\frac{3}{5}}{10}&=& \frac{\frac{3}{10}}{5}= \frac{\frac{1}{10}}{5}+\frac{\frac{1}{5}}{10}= \frac{1}{25}+\frac{1}{50},\text{ and}\\
    \frac{\frac{5}{7}}{8} &=&\frac{\frac{5}{8}}{7}=\frac{\frac{1}{2}}{7}+ \frac{\frac{1}{8}}{7}.
    \end{eqnarray*}
    Fibonacci points out that the composite fraction $\frac{\frac{5}{10}}{8}$ should not be expanded like $\frac{\frac{4+1}{8}}{10}$.
 \end{itemize}

\subsection*{The third category}
The third category is treating three different cases.
\begin{itemize}
\item The first case in third category is when $p$ divides $q+1$, i.e. there exists $k>1$ so that $kp=q+1$. Then
    \[\frac{p}{q}= \frac{1}{k}+ \frac{\frac{1}{k}}{q}.\]
    This derivation is also found in \cite[p.82]{Luneburg1993}. A more common \cite[p.13]{Cantor1900} \cite[p.620]{Sigler2002} way to write this is for $q=kp-1$
    \[\frac{p}{kp-1}=\frac{1}{k}+ \frac{\frac{1}{k}}{kp-1}.\]
    Fibonacci gives the following 5 examples.
    In the first example 2 divides 11+1 so
    \[\frac{2}{11} = \frac{1}{6}+\frac{\frac{1}{6}}{11}, \]
    and in the second example 3 divides 11+1
      \[\frac{3}{11} = \frac{1}{4}+\frac{\frac{1}{4}}{11} = \frac{1}{4}+\frac{1}{44}. \]
    In the next two examples 4 divides 11+1, and 6 divides 11+1 so
    \[\frac{4}{11} = \frac{1}{3}+\frac{\frac{1}{3}}{11} = \frac{1}{3}+\frac{1}{33}, \text{and }
      \frac{6}{11} = \frac{1}{2}+\frac{\frac{1}{2}}{11} = \frac{1}{2}+\frac{1}{22}.\]
    The final example is $\frac{5}{19}$ and 5 divides 19+1 so
    \[\frac{5}{19} = \frac{1}{4}+\frac{\frac{1}{4}}{19} = \frac{1}{4}+\frac{1}{76}.\]
\item The second case is a composite case $\frac{\frac{p}{q}}{r}$ and $\frac{p}{q}$ has a known unit fraction expansion. The two examples from Fibonacci are
      \[ \frac{\frac{2}{3}}{7}= \frac{\frac{1}{2}}{7}+ \frac{\frac{1}{6}}{7}  \]
      using the known expansion $\frac{2}{3}= \frac{1}{2}+\frac{1}{6}$ and
    \[ \frac{\frac{4}{7}}{9}= \frac{\frac{1}{2}}{9}+ \frac{\frac{1}{14}}{9}  \]
    using $\frac{4}{7}= \frac{1}{2}+\frac{1}{14}$.
 \item The next is a reverse composite case
    \[\frac{\frac{p}{q}}{r}=\frac{\frac{p}{r}}{q}.\]
 The two examples are
 \[\frac{\frac{3}{7}}{11}=\frac{\frac{3}{11}}{7}= \frac{\frac{1}{4}}{7}+\frac{\frac{1}{44}}{7},\]
 by the third rule (category) (3 divides 1+11) and
  \[\frac{\frac{3}{7}}{8}=\frac{\frac{3}{8}}{7}\]
 which belongs to two categories. For the second category $\frac{3}{8}$, 3=1+2 where each number in the sum divides 8 so  $\frac{3}{8}= \frac{2}{8}+\frac{1}{8}$, and for the third category 3 divides 8+1 so $\frac{3}{8}= \frac{1}{3}+\frac{\frac{1}{3}}{8}$.

\item  This paragraph is called {\em On the same category/distinction} \cite{Dunton1966,Sigler2002}.

If $q+1$ is a practical number then for $1<p<q$ the numerator  $p$ is a sum of divisors of $q+1$, using the second category. Then the fraction $\frac{p}{q}$ can be written as a sum  where each term  will have a two term unit fraction expansion using the third category. Let $s_i$ be divisor of $q+1$, then there exist $k_is_i=q+1$ and $p$ written as a sum of the divisors $p=\sum_{i=1}^{\ell}s_i$, then
\[\frac{p}{q}=\sum_{i=1}^{l}\frac{s_i}{q}=\sum_{i=1}^{\ell} \left(\frac{1}{k_i}+\frac{\frac{1}{k_i}}{q}\right).\]

For the three fractions $\frac{8}{11}$, $\frac{9}{11}$, and $\frac{10}{11}$ Fibonacci writes the numenator as 6+2, 6+3, and 6+4 respectively. Then
    \[\frac{8}{11}= \frac{6}{11}+\frac{2}{11}= \frac{1}{2}+\frac{1}{22}+\frac{1}{6}+\frac{1}{66},  \]
    \[\frac{9}{11}= \frac{6}{11}+\frac{3}{11}= \frac{1}{2}+\frac{1}{22}+\frac{1}{4}+\frac{1}{44},  \]
    and
    \[\frac{10}{11}= \frac{6}{11}+\frac{4}{11}=
      \frac{1}{2}+\frac{1}{22}+\frac{1}{3}+\frac{1}{33}.  \]

\end{itemize}

 \subsection*{The fourth category}
 In the fourth category $p-1$ divides $q+1$, then there exist $k>1$ so that
 $(p-1)k=q+1$. Here Fibonacci assumes  $q$ is prime.  Consider
 \[\frac{p}{q}= \frac{p-1}{q}+ \frac{1}{q}= \frac{1}{k}+ \frac{\frac{1}{k}}{q}+ \frac{1}{q}, \]
 from the third category \cite[p.83]{Luneburg1993}. An alternative  formulation \cite[p.13]{Cantor1900} is to let $p=p'+1$ then $p'k=q+1$
 \[\frac{p}{q}=\frac{p'+1}{p'k-1}= \frac{1}{k}+\frac{1}{p'k-1}+\frac{1}{k(p'k-1)}.\]
 Fibonacci demonstrates this rule with a subtraction step
  \[\frac{p}{q}-\frac{1}{q}= \frac{1}{k}+ \frac{\frac{1}{k}}{q}, \]
 with 5 examples. In the first example 4 divides 11+1, in the second  6 divides 11+1 and then apply the third rule. The other examples follows using the rule in third category .
 \begin{eqnarray*}
 \frac{5}{11} - \frac{1}{11}&=&\frac{4}{11}=
      \frac{1}{3}+\frac{\frac{1}{3}}{11},\text{  so  } \frac{5}{11} = \frac{1}{3}+\frac{1}{11}+\frac{1}{33}.\\
  \frac{7}{11} - \frac{1}{11} &=&\frac{6}{11}=
       \frac{1}{2}+\frac{\frac{1}{2}}{11},\text{  so  } \frac{7}{11}=\frac{1}{2}+\frac{1}{11}+\frac{1}{22}.\\
 \frac{3}{7} - \frac{1}{7}&=&\frac{2}{7}=
       \frac{1}{4}+\frac{\frac{1}{4}}{7},\text{  so  } \frac{3}{7} =\frac{1}{4}+\frac{1}{7}+\frac{1}{28}.\\
 \frac{6}{19}-\frac{1}{19}&=&\frac{5}{19}=
 \frac{1}{4}+\frac{\frac{1}{4}}{19}\text{ so } \frac{6}{19} =\frac{1}{4}+\frac{1}{19}+\frac{1}{76}\\
 \frac{7}{29} -\frac{1}{29}&=& \frac{6}{29}= \frac{1}{5}+\frac{\frac{1}{5}}{29},\text{  so  } \frac{7}{29}= \frac{1}{5}+\frac{1}{29}+\frac{1}{145}.
 \end{eqnarray*}

 \subsection*{The fifth category}
 In the fifth category $q$ is even and $p-2$ divides $q+1$. Then there exists $k>1$ so that $k(p-2)=q+1$ and
     \[\frac{p}{q}= \frac{2}{q} + \frac{p-2}{q}=\frac{1}{q/2}+\frac{1}{k}+\frac{\frac{1}{k}}{q}\]
     since the term $\frac{p-2}{q}$ belongs to the third category.
     Since $q$ is an even number the three numbers $q+1$, $k$, and $p-2$  are odd numbers. Let $p-2=2a+1$ and  $k=2n+1$ for positive integers $a$ and $n$. Then, using $q=k(p-2)-1$, we have the expression
     \begin{eqnarray*}
     \frac{p}{q}&=&\frac{2a+3}{(2n+1)(2a+1)-1}\\
     &=& \frac{1}{2n+1}+\frac{1}{(2n+1)a+n}+\frac{1}{(2n+1)[(2n+1)(2a+1)-1]}
     \end{eqnarray*}
     which gives the expression given by Cantor \cite[p.13]{Cantor1900}. Fibonacci demonstrates this rule using a subtraction
     \[\frac{p}{q}- \frac{2}{q} = \frac{p-2}{q}=\frac{1}{k}+\frac{\frac{1}{k}}{q}\]
     The two examples illustrating this rule are
     \begin{eqnarray*}
 \frac{11}{26}- \frac{2}{26}&=&\frac{9}{26}=
      \frac{1}{3}+\frac{\frac{1}{3}}{26},\text{ so } \frac{11}{26} =\frac{1}{3}+\frac{1}{13}+\frac{1}{78},\text{ and}\\
  \frac{11}{62} -\frac{2}{62}&=&\frac{9}{62}= \frac{1}{7}+\frac{\frac{1}{7}}{62}\text{ so } \frac{11}{62} =\frac{1}{7}+\frac{1}{31}+\frac{1}{434}.
      \end{eqnarray*}

\subsection*{The sixth category}
In the sixth category $q$ is divisible by 3 and $p-3$ divides $q+1$. Then there exists $k>1$ so that $k(p-3)=q+1$ and
     \[\frac{p}{q}= \frac{3}{q} + \frac{p-3}{q}=\frac{1}{q/3}+\frac{1}{k}+\frac{\frac{1}{k}}{q}.\]
Again, Fibonacci take a subtraction and then use category 3.     The two examples illustrating this rule are
          \begin{eqnarray*}
 \frac{17}{27} -\frac{3}{27}&=&\frac{14}{27}= \frac{1}{2}+\frac{\frac{1}{2}}{27}\text{ so } \frac{17}{27}=\frac{1}{2}+\frac{1}{9}+\frac{1}{54}, \text{ and}\\
  \frac{20}{33}-\frac{3}{33} &=&\frac{17}{33}=
     \frac{1}{2}+\frac{\frac{1}{2}}{33}\text{ so } \frac{20}{33}=\frac{1}{2}+\frac{1}{11}+\frac{1}{66}.
      \end{eqnarray*}

\subsection*{The seventh category}
The seventh Category is when none of the above rules or categories are applicable and general rules are described. The Category is divided in two parts, first the greedy method that will always work, but combined with the other categories will give neater results than just applying the greedy rule. The second part is what Fibonacci calls the universal rule.

\subsubsection*{The greedy algorithm}

      \begin{itemize}
      \item The first example Fibonacci considers is the rational $\frac{4}{13}$. The largest unit fraction that can be subtracted from  $\frac{4}{13}$  is $\frac{1}{4}$. The remainder is
          \[\frac{4}{13}-\frac{1}{4}=\frac{3}{4\cdot13}= \frac{\frac{3}{4}}{13}=\frac{\frac{1+2}{4}}{13}=\frac{1}{52}+\frac{1}{26},\]
           by the second category.
          But the remainder $\frac{\frac{3}{4}}{13}=\frac{3}{52}=\frac{1+2}{52}=\frac{1}{52}+\frac{1}{26}$ by the second category. The fraction
          \[\frac{4}{13}= \frac{1}{4}+\frac{1}{26}+\frac{1}{52}.\]
          However, for the fraction $\frac{3}{52}$ the largest unit fraction that can be subtracted is $\frac{1}{18}$ and the remainder is
          \[\frac{3}{52}-\frac{1}{18}= \frac{3-\frac{52}{18}}{52}=\frac{\frac{1}{9}}{52}=\frac{1}{468}.\]
          The unit fraction expansion of the fraction using the greedy algorithm is then
          \[\frac{4}{13}= \frac{1}{4}+\frac{1}{18}+\frac{1}{468}.\]

      \item For the next example $\frac{9}{61}$, Fibonacci gives two different expansions. First compute the largest unit fraction in $\frac{9}{61}$ (which is 1/7) and the consider the remaining part
          $\frac{9}{61}-\frac{1}{7}=\frac{9-\frac{61}{7}}{61}=\frac{2}{7\cdot61}=\frac{2}{427}$. Now using the third category
          \[\frac{2}{427}= \frac{1}{214}+ \frac{\frac{1}{214}}{427},\]
          Then $\frac{9}{61}=\frac{1}{7}+\frac{1}{214}+\frac{1}{214\cdot427}$.

          For the second way to do the expansion, consider the remaining part
          \(\frac{2}{427}=\frac{2}{7\cdot61}=\frac{\frac{2}{7}}{61}\) and using third composite
           category the well known \(\frac{2}{7}=\frac{1}{4}+\frac{1}{28}\).
           So
           \[\frac{9}{61}-\frac{1}{7}=\frac{\frac{2}{7}}{61}=
              \frac{\frac{1}{4}}{61}+\frac{\frac{1}{28}}{61}=\frac{1}{244}+\frac{1}{1708}.\]
           \item Third example is $\frac{17}{29}$. This is the only example where Fibonacci is treating the remainder using the greedy approach. Largest unit fraction is $\frac{1}{2}$ and the remainder is $\frac{17}{29}-\frac{1}{2}= \frac{2\frac{1}{2}}{29}= \frac{5}{58}$. Largest unit fraction in $\frac{5}{58}$ is $\frac{1}{12}$ and the remainder is \[\frac{5}{58}-\frac{1}{12}=\frac{5}{58}-\frac{\frac{58}{12}}{58}=
               \frac{5-4\frac{5}{6}}{58}=\frac{1}{6\cdot58}=\frac{1}{348}.\]
               The expansion is then
           \[\frac{17}{29}=\frac{1}{2}+\frac{1}{12}+\frac{1}{348}.\]
      \end{itemize}
\subsubsection*{A universal rule}
Sigler \cite[p.124]{Sigler2002} translates this other case as {\em a universal rule for separation into unit fractions}.
      Cantor \cite[p.13]{Cantor1900} states that the other case of seventh category is likely the most useful technique.

      Let $m$ be a large positive integer with many factors. Fibonacci exemplifies this by $m$ equal 12, 24, 36, 48 and 60 where $\frac{q}{2}<m$.
      He is also restricting the size of $m$, $m<2q$. Then
      \[\frac{p}{q}= \frac{m \;\cdot\;\frac{p}{q}}{m}= \frac{\lfloor \frac{mp}{q}\rfloor}{m} + \frac{\frac{s}{q}}{m},\]
      where $s$ is an integer $0\leq s <q$. If $r$ divides $m$ or can be written as a sum of divisors of $m$, then according to the reverse composite of the first or second Category, the second term is a unit fraction or a sum of unit fractions. Further if the integer $\lfloor \frac{mp}{q}\rfloor$ divides $m$ or can be written as a sum of divisors of $m$ then $\frac{p}{q}$ has a unit fraction expansion.
      
      In the first example Fibonacci consider the rational $\frac{17}{29}$ for $m=24$ and $m=36$ and using the second category for the term $\frac{\frac{s}{q}}{m}$.
      \begin{itemize}
      \item For the first example $m=24$ and $m=36$ and the fraction is $\frac{17}{29}$.
 \[\frac{17}{29}=\frac{\frac{24\cdot17}{29}}{24} = \frac{14\frac{2}{29}}{24}, \]
 and \(\frac{14}{24}\) is equal \(\frac{1}{3}+\frac{1}{4}\) or \(\frac{1}{2}+\frac{1}{12}\) using $14=6+8$ or $14=2+12$ and the second category. The second term is $\frac{\frac{2}{29}}{24}=\frac{\frac{2}{24}}{29}= \frac{1}{12\cdot29}=\frac{1}{348}$.
 So
 \[\frac{17}{29}=\frac{1}{3}+\frac{1}{4}+\frac{1}{348} \text{ or }
 \frac{17}{29}=\frac{1}{2}+\frac{1}{12}+\frac{1}{348}.\]

 For the case $m=36$
 \[\frac{ \frac{36\cdot17}{29}}{36} = \frac{21\frac{3}{29}}{36}. \]
From the integral part of the numerator \(\frac{21}{36}=\frac{7}{12}\) is equal \(\frac{1}{3}+\frac{1}{4}\) or equal
\(\frac{1}{2}+\frac{1}{12}\) (this follows from $7=6+1$ or $7=4+3$ and the second category.) For the fractional part \(\frac{\frac{3}{29}}{36}=\frac{\frac{3}{36}}{29}=\frac{\frac{1}{12}}{29}=\frac{1}{348}\)
which gives the same unit fraction expansion as for $m=24$.

\item The fraction $\frac{19}{53}$ is expanded in two ways. First using the fourth category
\[\frac{19}{53}=\frac{1}{3}+\frac{1}{53}+\frac{1}{159}.\]
Using the universal rule, first compute the largest unit fraction in $\frac{19}{53}$ (which is 1/3) and consider the remainder

\[\frac{19}{53} - \frac{1}{3}=\frac{19 -\frac{53}{3}}{53}=\frac{19 - 17\frac{2}{3}}{53}=
\frac{1\frac{1}{3}}{53}=\frac{1}{53}+\frac{1}{159}.\]

\item Next example is the fraction $\frac{20}{53}$. Multiply numenator and denominator with $m=48$
    \[\frac{20\cdot48}{53\cdot48}=\frac{18\frac{6}{53}}{48}=\frac{18}{48}+\frac{\frac{1}{8}}{53}.\]
    From the second category and the table of separation for $q=24$ we have
    \(\frac{18}{48}=\frac{9}{24}=\frac{1}{4}+\frac{1}{8} \text{ or } \frac{18}{48}=\frac{1}{3}+\frac{1}{24}.\)
\end{itemize}
\subsection*{Elegant and not so elegant}
The final part of Chapter 7 in Liber Abaci is  where Fibonacci discusses different unit fraction expansions for which some are more elegant than others. The first observation is to look at different fraction before the expansion into unit fractions, like if the rational is not reduced, like $\frac{6}{9}=\frac{2}{3}$ which belongs to the third category or
$\frac{6}{8}=\frac{3}{4}$ and using the technique in the second category. Further, simplify expressions like $\frac{3\frac{1}{2}}{8}=\frac{7}{16}$. Similarly,
\[\frac{4\frac{3\frac{2}{3}}{5}}{9}=\frac{71}{135}.\]
Largest unit fraction will be $\frac{1}{2}$ and the remainder is $\frac{71}{135}-\frac{1}{3}=\frac{7}{270}= \frac{1+6}{270}= \frac{1}{45}+\frac{1}{270}$.
Hence,
\[\frac{71}{135}=\frac{1}{2}+\frac{1}{45}+\frac{1}{270}.\]

Fibonacci points out that when the largest unit fraction is used in the general rule, the expansion of the remainder may be less than elegant. To make it more elegant, choose not the smallest denominator but one that is one or two larger. So for the first unit fraction of $\frac{p}{q}$   choose
\[\lceil \frac{q}{p}\rceil+1 \text{ or } \lceil \frac{q}{p}\rceil+2.\]
Consider the fraction $\frac{4}{49} $, then $\lceil \frac{49}{4}\rceil=13$ and
\( \frac{4}{49}-\frac{1}{13}=\frac{3}{637}=\frac{1}{319}+\frac{1}{637}+\frac{1}{319\cdot637}
\) by the fourth category. A result that is less than elegant. Consider instead
\( \frac{4}{49}-\frac{1}{14}=\frac{1}{98}\) which is more elegant.

An alternative derivation is given by
\[ \frac{4}{49}=\frac{\frac{4}{7}}{7}=\frac{\frac{1}{2}+\frac{1}{14}}{7}= \frac{1}{14}+\frac{1}{98},\]
by the third category.

\subsection*{Table of separations}\label{sec:separation}
Let $q$ be a practical number then for every $p$, $1<p<q$, there exist divisors of $q$ so that $p$ is a sum of distinct divisors. Let \(p=\sum_{i=1}^{s} q_i\) where $q_i$ is a divisor of $q$, ie for some integer $d_i$ we have $d_iq_i=q$, and the fraction $\frac{p}{q}$ is an Egyptian fraction
\begin{equation}\label{eq:sep}
  \frac{p}{q}= \sum_{i=1}^{s}\frac{1}{d_i}.
\end{equation}
Consider the practical number $q=12$ where the divisors are $1,2,3,4,6,12$. For $p=5$ then
$p$ can be written as a sum of two divisors in different ways $p=2+3$ or $p=1+4 $ which gives
\[\frac{5}{12}=\frac{1}{6}+\frac{1}{4}=  \frac{1}{3}+\frac{1}{12}.\]
Fibonacci is choosing the first of these. In the tables of separation $q$ is $6,8,12,20,24,60,$ and $100$ and Fibonacci gives all $1<p<q$ for $q<60$ and only some for $q=60$ and $q=100$. In the tables in this paper we only consider the cases where $\frac{p}{q}$ can not be reduced to a unit fraction. For completeness, the divisors are given
\begin{eqnarray*}
  q &=6& \text{ the divisors are }  1,2,3,6 \\
  q &=8& \text{ the divisors are }  1,2,4,8\\
  q &=12& \text{ the divisors are }  1,2,3,4,6,12\\
  q &=20& \text{ the divisors are }  1,2,4,5,10,20\\
  q &=24& \text{ the divisors are }  1,2,3,4,6,8,12,24\\
  q &=60& \text{ the divisors are }  1,2,3,4,5,6,10,12,15,20,30,60\\
  q &=100& \text{ the divisors are }  1,2,4,5,10,20,25,50,100.\\
\end{eqnarray*}
\begin{center}
\begin{tabular}{|l||l|}
  \hline
  $q=6$&         $q=12$\\ \hline
  $p=4=3+1$ &    $p=5=3+2=4+1$ \\
  $p=5=3+2$ &    $p=7=6+1$ \\\cline{1-1}
  $q=8$     &    $p=8=6+2$\\\cline{1-1}
  $p=3=2+1$ &    $p=9=6+3$ \\
  $p=5=4+1$ &    $p=10=6+4$ \\
  $p=6=4+2$ &    $p=11=6+3+2=6+4+1$ \\\cline{2-2}
  $p=7=4+2+1$&      \\\cline{1-1}

\end{tabular}\end{center}
\begin{center}
  \begin{tabular}{|l|l|}
    \hline
    \multicolumn{2}{|c|}{$q=20$}\\\hline
    $p=3=2+1$       &$p=13=10+2+1$\\
    $p=6=5+1=4+2$   &$p=14=10+4$\\
    $p=7=5+2$       &$p=15=10+5$\\
    $p^\dag=8$      &$p=16=10+4+2=10+5+1$\\
    $p=9=5+4$       &$p=17=10+5+2$\\
    $p=11=10+1$     &$p^\dag=18$\\
    $p=12=10+2$     &$p=19=10+5+4$\\

    \hline
  \end{tabular}
\end{center}
Notes:
\begin{itemize}
\item For $p=8$ (and $q=20$) 8 is not a sum of two of the divisors of 20. However, 8 is a sum of three distinct divisors $8=5+2+1$. Fibonacci writes $\frac{2}{5}$. Choose $m=3$ and consider
$\frac{p\cdot m}{q\cdot m}$, and $p\cdot m=24=20+4$ which is found in the table $q=60$.
\item For $p=18$ (and $q=20$) 18 is not a sum of two or three of the divisors of 20. However, 18 is a sum of four distinct divisors $18=10+5+2+1$. Fibonacci writes $\frac{18}{20}=\frac{1}{15}+\frac{1}{3}+\frac{1}{2}$.
Choose $m=3$ and consider $\frac{p\cdot m}{q\cdot m}$. Then $p\cdot m= 30+20+4$ and
$\frac{54}{60}= \frac{1}{15}+\frac{1}{3}+\frac{1}{2}$. An entry for $p=54$ is not included in the table for $q=60$.
\end{itemize}
\begin{center}
  \begin{tabular}{|l|l|}
    \hline
    \multicolumn{2}{|c|}{$q=24$}\\\hline
    $p=5=3+2=4+1$       &$p=16=12+4$\\
    $p=7=4+3=6+1$       &$p=17=12+3+2=$\\
                        &\quad$=12+4+1=8+6+3$\\
    $p=9=6+3=8+1$       &$p=18=12+6$\\
    $p=10=6+4=8+2$      &$p=19=12+4+3=12+6+1$\\
    $p=11=8+3$       &$p=20=12+8$\\
    $p^{\dag\dag}=13=6+4+3=$     &$p=21=12+6+3=12+8+1$\\
    \quad\quad$=8+4+1=8+3+2$&\\
    $p=14=8+6=12+2$     &$p=22=12+6+4=12+8+2$\\
    $p=15=12+3$         &$p=23=12+8+3$\\
    \hline
  \end{tabular}
\end{center}
Note:
\begin{itemize}
\item For $p=13$ (and $q=24$) Fibonacci writes $13=6+4+3$. However,   $p=12+1$ is a two term expansion.
\end{itemize}
For $q=60$, Fibonacci consider  $1\leq p\leq 31$, and $p=35,40,50,55$. In the revised version of Liber Abacci $p=55$ is not included.
\begin{center}
  \begin{tabular}{|l|l|}
    \hline
    \multicolumn{2}{|c|}{$q=60$}\\\hline
    $p=7=6+1=5+2=4+3$       &$p=23=20+3$\\
    $p=8=6+2=5+3   $        &$p=24=20+4$\\
    $p=9=6+3=5+4$           &$p=25=20+5=15+10$\\
    $p=11=10+1=6+5$         &$p=26=20+6$\\
    $p=13=10+3=12+1$        &$p=27=15+12$\\
    $p=14=10+4=12+2$        &$p=28=12+10+6=20+6+2$\\
    $p=16=10+6=12+4=15+1$   &\ \ $=20+5+3=15+12+1$\\
                            &\ \ $=15+10+3$\\
    $p=17=15+2=12+5$        &$p=29=20+6+3=20+5+4$\\
                            &\ \ $=15+12+2=15+10+3$\\
    $p=18=12+6=15+3$        &$p=31=30+1$\\
    $p=19=15+4$             &$p=35=20+15=30+5$\\
    $p=21=15+6=20+1$        &$p=40=30+10$\\
    $p=22=20+2=12+10$       &$p=50=30+20$\\
                            &$p=55=30+15+10=30+20+5$\\
    \hline
  \end{tabular}
\end{center}
\begin{center}
  \begin{tabular}{|l|l|}
    \hline
    \multicolumn{2}{|c|}{$q=100$}\\\hline
    $p=3=2+1$               &$p^{\dag\dag\dag}=60$\\
    $p=6=4+2=5+1   $        &$p=70=20+3$\\
    $p=7=5+2=5+4$           &$p^{\dag\dag\dag}=75$\\
    $p^{\dag\dag\dag}=8$    &$p^{\dag\dag\dag}=80$\\
    $p=9=5+4$               &$p^{\dag\dag\dag}=85=50+25+10$\\
    $p=15=10+5$             &$p=95=50+25+20$\\
    $p=30=20+10$            &$p=96=50+25+20+1$\\
    $p=35=25+10$            &$p=97=50+25+20+2$\\
    $p^{\dag\dag\dag}=40$   &$p=98=50+25+20+2+1$\\
    $p=45=25+20$            &$p=99=50+25+20+4$\\
    \hline
  \end{tabular}
\end{center}
Notes on lines marked $p^{\dag\dag\dag}$ in the table for $q=100$.
\begin{itemize}
\item For $p=8$ (and $q=100$) Fibonacci writes $\frac{2}{25}$. There are no two terms expansion of 8 in terms of the divisors of 100. However, there are two terms $\frac{2}{25}=\frac{1}{15}+\frac{1}{75}$. A possible three terms expansion could be based on $8=5+2+1$
\item For $p=40$ (and $q=100$) Fibonacci writes $\frac{2}{5}$. Since $\frac{40}{100}=\frac{2}{5}=\frac{24}{60}$ the expansion is found in the table for $q=60$ in the line for $p=24$.
 \item For $p=60$ (and $q=100$) Fibonacci writes $\frac{3}{5}$, but $\frac{3}{5} = \frac{12}{20}=\frac{1}{2}+\frac{1}{10}$ in the table for $q=20$ line $p=12$. A possible expansion could be $60=50+10$.
\item For $p=75$ (and $q=100$) Fibonacci writes $\frac{3}{4}$. The expansion $\frac{3}{4}=\frac{6}{8}=\frac{1}{2}+\frac{1}{4}$ can be found in the table for $q=8$ and line $p=6$.
    \item For $p=80$ (and $q=100$) Fibonacci writes $\frac{4}{5}$. But $\frac{4}{5}=\frac{16}{20}=\frac{1}{2}+\frac{1}{5}+\frac{1}{10} $ from table for $q=20$ line $p=16$. Possible expansions in terms of the divisors of $q=100$ could be $80=50+20+10=50+25+5$.
\item For $p=85$ (and $q=100$) Fibonacci writes $\frac{3}{4}+\frac{1}{10}$. An alternative approach would be
$\frac{85}{100}=\frac{17}{20}=\frac{1}{2}+\frac{1}{4}+\frac{1}{10} $ from the table $q=20$ and line $p=17$. In the revised version the expansion is based on $85=50+25+10$.
\end{itemize}

Fibonacci is not consistently choosing the shortest possible unit fraction expansion or choosing the fraction with smallest largest denominator.  A common belief is that the last term in the expansion should have as small denominator as possible or the shortest possible expansion. However, as stated, the tables of separation does not support this.

\section*{Concluding remarks}
This paper is based on the a web-page \href{https://liberabaci.blogspot.com/}{{\em Liber Abaci and seven Egyptian fraction methods}} by Milo Gardner January 4, 2007.

\bibliographystyle{plain}

\end{document}